\documentclass[a4paper]{article}
\title{Successful strategies for a queens placing game on
 an $n$ x $n$ chess board}
\author{Thomas Jenrich}
\date{2014-04-21}
\begin{document}
\maketitle

\section{Abstract and introduction}

In his list of open problems (\cite{Eri}), Martin Erickson described
a certain game:

``Two players alternately put queens on an $n$ x $n$ chess board so that each new queen is not in range of any queen
  already on the board (the color of the queens is unimportant). The last player who can move wins.''

Then he asked: ``Who should win?''

Obviously, for $n$ up to 3, the first player wins, if he
does not miss to start at the central position in the case $n=3$.

In this article, we give very simple always winning strategies for the first
player if $n$ is 4 or odd.
The additionally (in the source package) provided computer program QPGAME3 has
been used to check that there are successful strategies for the first player
if $n$ is 6 or 8, and for the second player if $n$ is 10, 12, 14, or 16.

As discovered during the submission process of the first version of this
article, Hassan A Noon presented consistent results concerning values of
$n$ which are odd or at most 10, in his B.A. thesis (\cite{Noo1}) and,
together with Glen Van Brummelen, in a journal article (\cite{Noo2}).

\section{The case $n=4$}

The first player places a queen at one of the four inner positions.
When * indicates an unavailable cell and . indicates an available cell,
we do have this situation (up to meaningless rotations of the complete board in
steps of 90 degrees):

\ttfamily

 * * * .

 * * * *

 * * * .

 . * . *

\rmfamily

No matter which of the four available cells the second player takes, it
remains exactly one available cell for the first player to place the
last queen.

\section{Notation and conflict-free positions}

The rows and columns of the board are numbered from 0 to $n-1$. A position
is an (ordered) pair of a row number $r$ and a column number $c$, written
as ($r$,$c$).

If two queens are placed at ($r_1$,$c_1$) and ($r_2$,$c_2$), then the
game rules require:

 $r_1 \ne r_2$ (different rows),

 $c_1 \ne c_2$ (different columns),

 $r_1+c_1 \ne r_2+c_2$ (different falling diagonals), and

 $r_1-c_1 \ne r_2-c_2$ (different rising diagonals).

\section{The case of odd $n$}

The first player puts the first queen to the cell at the central position
($(n-1)/2$,$(n-1)/2$).

From then on: When the second player has put a queen to the cell at a position
($r$,$c$), the first player puts a queen to the cell at position ($R$,$C$),
where $R=n-1-r$ and $C=n-1-c$. This is a valid position, received from
($r$,$c$) by mirroring the two coordinate values with respect to the middle
lines, or, with the same effect, by rotating the board a half turn.

This strategy is always successful for the following reasons:

The central cell lies just at the middle of the straight line between
($r$,$c$) and ($R$,$C$). Thus, a conflict between queens at these positions
exists if and only if there is a conflict between ($r$,$c$) and the queen at
the central position; but this is forbidden by the game rules.
If ($R$,$C$) was available before a queen is put to ($r$,$c$), then it is
still available after that move.

After the placing of the first queen, there is (half turn) rotational
symmetry. Because of this symmetry, ($R$,$C$) is available if and only
if ($r$,$c$) is available. After putting two queens to the cells at these
positions, the (half turn) rotational symmetry exists again.

Eventually, no cell is available (for the second player) and the first player wins.

\section{The provided computer program QPGAME3}

As like as its two predecessors, the source code file QPGAME3.PAS has been
developed for PASCAL compilers compatible with Turbo Pascal 4.0. Lines are at
most 78 characters long. For inspections the use of an ASCII compatible
monospaced font is strongly recommended. The intended indentation is by one
character per structure level, using blanks (instead of tabs).

The program includes a good portion of comments (enclosed in curly braces).
So it should be fairly understandable at least by readers knowing at least
one imperative programming language.

Basically, the core procedure \emph{wins} performs a general backtracking
search. Of course, all selected positions have to be legal.

When a player in a certain situation finds a winning move position $P_1$ as
an answer to his opponent's move to a position $P_0$, it is clear that $P_0$
would be a winning answer to $P_1$ if the situation before both moves did not
change. Therefore such a winning move position $P_1$ is registered in the
internal variable \emph{forbidden} and skipped in further checks. That would
give no advantage if $P_1$ has been checked before $P_0$. Therefore,
excluding (``for historical reasons'') the second move of a game, the
search of an answering move starts in a row that follows the row of the move
of the previous move with respect to the search order for the previous move.

The newly added Boolean routine parameter \emph{rotsym} indicates, whether
all currently involved moves of player 2 had re-established the initial
(half-turn) rotation symmetry by choosing the unique position just a
half-turn away from the position of the previous move of player 1.
If the value of \emph{rotsym} is \emph{True} when player 1 has to make a
move, the symmetry is used by just considering the rows up to the center of
the board.
(Source code fragment: ``\verb+(not rotsym or (r<=nm1div2))+'').

In addition, there are these further conditions (in order to reduce
the set of positions to be checked):

For the first move, the (then complete) symmetry is used.

(Source code fragment: ``\verb+(c<=nm1div2) and (r<=c)+'').

If $n$ is odd, the above described strategy for the first player is used.

If $n$ is even and not greater than 8, the first player is forced to start at
a certain of the four inner positions.
(Source code fragment: ``(r=nm1div2) and (c=nm1div2)'').

Tests have shown, that this restriction does not prevent the first player
from winning (in that case) while significantly reducing the sum of moves
performed during the search in the case $n=8$.

Finally, there is a special treatment of the case $n=16$ by player 2,
described in the following subsection. (Source code fragment:

 ``\verb# ((n<>16) or#

   \verb#   (n_moves>3) or#

   \verb#   (n_moves=3) and n16_player_2_round_2(prev_r,prev_c,r,c) or#

   \verb#   (n_moves=1) and n16_player_2_round_1(prev_r,prev_c,r,c))#'').

\subsection{Special treatment of the case $n=16$}

In order to reduce the otherwise very long execution time, QPGAME3.PAS uses
previously found and checked moves of player 2 to treat the case $n=16$ in a
special way: Two routines (\emph{n16\_player\_2\_round\_1} and
\emph{n16\_player\_2\_round\_2}) select unique moves of player 2 when he
has to answer to the moves of player 1 in the first two rounds.
\vspace{0.1in}

In \emph{n16\_player\_2\_round\_1}, the table \emph{T} contains for each
possible first move of player 1 a hexadecimally noted byte value, where \$FF
indicates an invalid (unexpected) move of player 1 and each other value
encodes an answering move of player 2: The upper/lower half-byte
(first/second hexadecimal) gives the row/column number.

Because the first move of player 1 is restricted to positions satisfying
``\verb+(c<=nm1div2) and (r<=c)+'', here ``\verb+(c<=7) and (r<=c)+'', just
36 move positions have to bo considered. In order to be able to use the
rotation symmetry, it has been tried to use (15-r0,15-c0) as often as
possible. Clearly, this was impossible for the 8 positions on the diagonal
(r0=c0). In these cases a position (r1,c1) is given such that
``\verb+(c1<=7) and (r1<=c1)+'' and therefore two positions can and will be
checked by one actual check. So, 28 actual checks have to be executed.

By the way, for some non-diagonal positions (r0,c0) turned out that
player 2 could not win by choosing (15-r0,15-c0), for instance (1,2) and
(6,7).
\vspace{0.1in}

In \emph{n16\_player\_2\_round\_2}, the situation is more complicated, because
the answering second move of player 2 depends not just on the second move of
player 1 (third move of the game), but also on the first move of player 1
(first move of the game). Therefore, that first move has been stored before
in the global variables \emph{r0} and \emph{c0} and is used to address an
index value in the table \emph{A}. Here, the value 0 means, that the
addressing position is invalid (by design); its use would cause a range
check error and abort the execution. Otherwise, that value is the first index
value to address an entry in the table \emph{B}. As the other two index
values, the row and the column number of the second move of player 1
(routine parameters (\emph{r2} and \emph{c2}) are used. The type,
notation, meaning, and usage of the entries of \emph{B} are as described
above for the table \emph{T}.

\vspace{0.1in}

Some of the (sub-)tables contain unused space or redundant information. In
addition, the usage of the tables in a predicative (checking) instead of a
functional (constructive) way is not optimal for minimal runtime. But the
latter can be ignored because the selection of a second or fourth move of a
game is an extremely seldom event under the move selections performed during
a whole case check (search) if $n$ is not very little. And because the
program is still small, also the waste of space is no real problem here.

\subsection{General properties of QPGAME3}

The program does not use the heap or any pointer operation at all.
If you don't change the respective compiler directives, range checks
and stack overflow checks will be generated. So the resulting executable
will be extremely safe. It is also small and needs only a few
kilobytes for the stack.

The program ignores any command line parameters or inputs other than
pressing Ctrl-C to cancel the execution - where the speed of response
depends on the used compiler (slow in the case of Turbo Pascal).
It writes only to the standard output device and into the automatically
created or rewritten, resp., file QPGAME.LST within the working directory.
In the default case, the standard output device will be the monitor screen.
But one could redirect that output (e.g. to a file). In most cases that will
be unnecessary because the listing in QPGAME.LST is an essence of the data
written to the standard output device. Just some execution state indicators
are omitted. Even after the execution has been stopped by cancelling, the
listing file should be readable and contain the data written before the end
of the execution. During the execution, its content is not accessible (by
other processes).

\subsection{Progress indication and summaries}

As for the older program versions, a + will be emitted (written to the
standard output device) after each 1000000 moves. Those characters (and
corresponding new-lines) will not be written into the listing file.

When the check of a certain case (value of $n$) begins or ends, an
information line will be written (to the standard output decice and into the
listing file).

\subsubsection{The sum of calls}
QPGAME1, the program provided with the first version of this article, used
a simple longint (signed 32 bit integer) variable to count the calls (of the
core routine / of a player to try to find a winning move), then displayed
as \emph{Sum of moves}, which is actually smaller (by one). Because the upper
limit of that counter was $2^{31}-1 = 2147483647$, overflows (reductions by
$2^{32}$) could occur in realistic cases. Those overflows can not influence
the search but mislead the user. So, the counter has been redesigned in
QPGAME2.PAS. Since then, it will not overflow below
$(2^{31}-1) \times 10^6 > 2.14 \times 10^{15}$.

\subsubsection{Detailed sub-case information}

In order to allow to give more detailed information about the state and
subresults of the execution, especially in very extensive checks, now the
generation of additional output can be advised:

If the Boolean variable \emph{first\_move\_checking\_statistics} does have
the value \emph{True}, a subresult summary line (containing the position of
the first move (by player 1), the resulting outcome, and the sum of calls
during the check of that first move) will be written to the standard
output decice and into the listing file.

If the Boolean variable \emph{indicate\_third\_moves\_checking} does have
the value \emph{True}, the start and the end of the check of the situation
established by a third move will be indicated (but not written to the
listing file). To be more concrete: On the start of that check, the positions
of the first and third move of the game (both by player 1) will be written
to the standard output device. On the end of that check, the result of
the routine \emph{wins} (the sub-case success of player 2) will be indicated
by a digit (0 means loss, 1 means win). Typically, there are some + symbols
between these two outputs.

An example (from the check of the case $n=16$):
\vspace{0.1in}

\ttfamily
\small
[1:(0,0)]  3:(2,3)++++++++++++++++++++++++++++++++++++ -> 1
\rmfamily

\vspace{0.1in}

In the unchanged QPGAME3.PAS, ``\verb+n>=16+'' is the evaluating expression
for both mentioned control variables.

\subsection{Compiling and running, actually checked cases}

In principle, the program is able to check the cases of $n$ from 1 to 16 or
32 (depending on the compiler symbol \emph{BIGN}, see below) without large
changes.
Because some mentioned cases were already solved mathematically and
checking the case $n=14$ will take some time, in the published version only
checks for even values of $n$ from 6 to 12 are called from the main loop.

\vspace{0.1in}

In order to avoid a compilation result depending on the settings you could
use the command line versions of the compilers (TPC for Turbo Pascal, BPC for
Borland Pascal 7, DCC32 for Borland Delphi (32 bit versions; do not miss to
use the -CC option in order to generate a console executable), VPC for
Virtual Pascal, FPC for Free Pascal) instead of the compilers integrated
in the IDEs.

In order to make the execution faster than that of QPGAME1 (provided with
the first version of this article), alternative structure and usage of the
variables containing the information on already used rows, columns, and
diagonals as described in \cite{Ric} have been implemented with QPGAME2.
The program became indeed faster, but because the new code will not work if
$n$ is greater than 16 and its function is not just that obvious as that of
the old code, the old code is still in the source file and can be used instead
of the new code by compiling with the symbol BIGN defined. For this purpose,
the appropriate command line option would be -DBIGN for compilers from Borland
(Turbo Pascal, Delphi), and -dBIGN for Free Pascal.
\vspace{0.1in}

The program has been successfully compiled and executed on a 1 GHz Intel
PIII PC running MS Windows 98 SE. These are the used compilers and the
respective two execution times from compilations with defined/undefined
compiler symbol BIGN:

Turbo Pascal 5.5 : 1:44 min / 1:16 min

Turbo Pascal 7.01 : 1:39 min / 0:41 min

Borland Delphi 4.0 build 5.37 : 0:13 min / 0:11 min

Virtual Pascal 2.1 build 279 : 0:18 min / 0:14 min

Free Pascal 2.4.4 i386-Win32 : 0:15 min / 0:11 min

\vspace{0.1in}

Here is the content of QPGAME.LST after compiling the (unchanged)
QPGAME3.PAS and running the generated QPGAME3.EXE:

\vspace{0.1in}

\ttfamily
\small

=== Checking solutions for the queens placing game problem ===

=== \ \ Version 3 \ \ \ Copyright (c) 2014-04-17 Thomas Jenrich ===
\vspace{0.1in}

Hints:

\ Output listing into file QPGAME.LST within the working directory.

\ After each 1000000 moves a + will be emitted.

\ To cancel the execution press Ctrl-C.
\vspace{0.1in}

Starting search with n = 6

Search completed. Result of player 1: win. Sum of calls: 54
\vspace{0.1in}

Starting search with n = 8

Search completed. Result of player 1: win. Sum of calls: 2266
\vspace{0.1in}

Starting search with n = 10

Search completed. Result of player 1: loss. Sum of calls: 653007
\vspace{0.1in}

Starting search with n = 12

Search completed. Result of player 1: loss. Sum of calls: 11334613
\vspace{0.1in}

== Regular program stop ==

\rmfamily

\vspace{0.1in}

After changing the upper limit of the main loop from 13 to 14,
Free Pascal has compiled QPGAME3.PAS again. The execution of the
resulting QPGAME3.EXE took 18:51 min. This is the (relevant part of the)
additional output:
\vspace{0.1in}

\ttfamily
\small

Starting search with n = 14

Search completed. Result of player 1: loss. Sum of calls: 1161385667

\rmfamily
\vspace{0.1in}

After changing the lower limit and the upper limit of the main loop
into 16, Free Pascal has compiled QPGAME3.PAS again. The execution of the
resulting QPGAME3.EXE took 22:55:30 h. Parallel activities on that computer
may have caused some delay, but probably not more than one hour.
\vspace{0.1in}

Here is the content of the generated QPGAME.LST:
\vspace{0.1in}

\ttfamily
\small

=== Checking solutions for the queens placing game problem ===

===  \ \ Version 3 \ \ \ Copyright (c) 2014-04-17 Thomas Jenrich ===

\vspace{0.1in}

Hints:

\ Output listing into file QPGAME.LST within the working directory.

\ After each 1000000 moves a + will be emitted.

\ To cancel the execution press Ctrl-C.
\vspace{0.1in}

Starting search with n = 16

pl. 1: (0,0) -> pl. 2: win. Sum of calls: 4470810024

pl. 1: (0,1) -> pl. 2: win. Sum of calls: 3905839444

pl. 1: (0,2) -> pl. 2: win. Sum of calls: 3734401972

pl. 1: (0,3) -> pl. 2: win. Sum of calls: 3184149898

pl. 1: (0,4) -> pl. 2: win. Sum of calls: 4328139348

pl. 1: (0,5) -> pl. 2: win. Sum of calls: 4425220446

pl. 1: (0,6) -> pl. 2: win. Sum of calls: 2278166076

pl. 1: (0,7) -> pl. 2: win. Sum of calls: 2676717214

pl. 1: (1,2) -> pl. 2: win. Sum of calls: 3876758648

pl. 1: (1,3) -> pl. 2: win. Sum of calls: 2063718676

pl. 1: (1,4) -> pl. 2: win. Sum of calls: 2120404822

pl. 1: (1,5) -> pl. 2: win. Sum of calls: 2114679020

pl. 1: (1,7) -> pl. 2: win. Sum of calls: 2290066796

pl. 1: (2,2) -> pl. 2: win. Sum of calls: 3360372268

pl. 1: (2,3) -> pl. 2: win. Sum of calls: 1952040408

pl. 1: (2,4) -> pl. 2: win. Sum of calls: 1993953094

pl. 1: (2,5) -> pl. 2: win. Sum of calls: 1965208480

pl. 1: (2,6) -> pl. 2: win. Sum of calls: 2088818438

pl. 1: (2,7) -> pl. 2: win. Sum of calls: 2192529724

pl. 1: (3,4) -> pl. 2: win. Sum of calls: 1719300574

pl. 1: (3,5) -> pl. 2: win. Sum of calls: 1834911490

pl. 1: (3,6) -> pl. 2: win. Sum of calls: 1981765458

pl. 1: (3,7) -> pl. 2: win. Sum of calls: 1964816772

pl. 1: (4,5) -> pl. 2: win. Sum of calls: 1671384528

pl. 1: (4,6) -> pl. 2: win. Sum of calls: 1857679532

pl. 1: (4,7) -> pl. 2: win. Sum of calls: 2007652010

pl. 1: (5,6) -> pl. 2: win. Sum of calls: 1690434580

pl. 1: (5,7) -> pl. 2: win. Sum of calls: 1712035496

Search completed. Result of player 1: loss. Sum of calls: 71461975237
\vspace{0.1in}

== Regular program stop ==

\rmfamily

\vspace{0.1in}

Author's eMail address: thomas.jenrich@gmx.de (Deutsch, English)

\end{document}